\documentclass[letterpaper, 10pt,psfig, conference]{ieeeconf}

\usepackage{amssymb}
\usepackage{graphicx}
\usepackage{amsmath}%
\usepackage{amsfonts}
\usepackage{epsfig}

\newcommand{\nc}{\newcommand}
\nc{\dd}{d_{\Delta}}
\nc{\uk}{\underline{k}}
\nc{\bk}{\bar{k}}

\nc{\kib}{\bar{k}_i}
\nc{\kioneb}{\underline{k}_{i+1}}

\nc{\that}{\hat{\theta}}

\nc{\rhod}{{\rho_{\delta}}}
\nc{\rhot}{{\rho_{\delta} ( \phi (t-d+1) , e(t+1))}}
\nc{\rhoj}{{\rho_{\delta} ( \phi (j-d+1) , e(j+1))}}
\nc{\rhotm}{{\rho_{\delta} ( \phi (t-d+1) , e(t+1))}}
\nc{\rhojm}{{\rho_{\delta} ( \phi (j-d+1) , e(j+1))}}

\nc{\bphi}{\bar{\phi}}
\nc{\tphi}{\tilde{\phi}}
\nc{\tct}{{\tilde{\check{\theta}}}}
\nc{\vc}{\check{V}}
\nc{\ttt}{\tilde{\theta}}
\nc{\D}{{{\bf D}}}

\nc{\xnun}{{\tbo{x(0)}{\hat{u}^0 (0)}}}
\nc{\xkuk}{{\tbo{x(kT)}{\hat{u}^0(kT)}}}

\nc{\blz}{\bar{\lambda}_0}
\nc{\bgz}{\bar{\gamma}_0}
\nc{\bh}{\bar{h}}
\nc{\bn}{\bar{n}}
\nc{\barki}{{\bar{k}_i}}
\nc{\barkio}{{\bar{k}_{i-1}}}

\nc{\bub}{{\bar{\underline{B}}}}
\nc{\bcub}{{\bar{\underline{C}}}}

\nc{\bb}{{\bar{B}}}
\nc{\bc}{{\bar{C}}}
\nc{\ub}{{\underline{B}}}
\nc{\uc}{{\underline{C}}}
\nc{\ucb}{{\underline{\cal B}}}
\nc{\ucc}{{\underline{\cal C}}}
\nc{\htxe}{{\hat{\tilde{x}}^{\eps}}}
\nc{\htue}{{\hat{\tilde{u}}^{\eps}}}
\nc{\htye}{{\hat{\tilde{y}}^{\eps}}}
\nc{\rmr}{{\R^{\tilde{m} \times \tilde{r}}}}
\nc{\kde}{{ \tilde{K}_{\delta}^{\eps}}}
\nc{\kdep}{{ \tilde{K}_{\delta}^{\eps}(p)}}
\nc{\omc}{{{\cal O}_{m}(C,A)}}
\nc{\tkt}{{ [kT , (k+1)T)}}
\nc{\mbl}{{\frac{\bar{m}}{\ell}}}
\nc{\calo}{{\cal O}}
\nc{\calc}{{\cal C}}
\nc{\acl}{{A_{cl}}}
\nc{\acla}{{A_{cl(A,B,C)}}}
\nc{\aclep}{\tilde{A}_{cl}^{\eps} (p)}
\nc{\gabc}{{\gamma_{(A,B,C)}}}
\nc{\gabcb}{{\gamma_{( \bar{A} , \bar{B} , \bar{C} )}}}
\nc{\gabcbt}{{\tilde{\gamma}_{(\bar{A} , \bar{B} , \bar{C})}}}
\nc{\dabc}{{\delta_{(A,B,C)}}}
\nc{\dabcb}{{\delta_{( \bar{A} , \bar{B} , \bar{C} )}}}
\nc{\labc}{{\lambda_{(A,B,C)}}}
\nc{\labcb}{{\lambda_{(\bar{A}, \bar{B}, \bar{C})}}}
\nc{\acleb}{{A_{cl( \bar{A} , \bar{B} , \bar{C} )}^{\eps}}}
\nc{\abcp}{{(A,B,C) \in {\cal P}}}
\nc{\abcbp}{{(\bar{A}, \bar{B}, \bar{C}) \in {\cal P}}}
\nc{\hxe}{{\hat{x}^{\eps}}}
\nc{\hue}{{\hat{u}^{\eps}}}

\nc{\eabcb}{{\varepsilon_{(\bar{A}, \bar{B}, \bar{C})}}}
\nc{\feabc}{{F_{(A,B,C)}^{\varepsilon}}}
\nc{\fabce}{{F_{(A,B,C)}^{\varepsilon}}}
\nc{\fabcb}{{F_{(\bar{A}, \bar{B}, \bar{C})}}}
\nc{\fabcbe}{{F_{(\bar{A}, \bar{B}, \bar{C})}^{\varepsilon}}}
\nc{\fabc}{{F_{(A,B,C)}}}
\nc{\oh}{{\cal O}}
\nc{\xbb}{{\bar{\bar{x}}}}
\nc{\bp}{{\bar{p}}}
\nc{\fguy}{{\cal F}}
\nc{\eps}{{\varepsilon}}
\nc{\feps}{{\hat{f}_{\eps}}}
\nc{\xbeps}{{\bar{x}^{\eps}}}
\nc{\xheps}{{\hat{\bar{x}}^{\eps}}}
\nc{\xhe}{{\hat{{x}}^{\eps}}}
\nc{\xte}{{\tilde{{x}}^{\eps}}}
\nc{\uhe}{{\hat{{u}}^{\eps}}}
\nc{\ute}{{\tilde{{u}}^{\eps}}}

\nc{\nabcb}{{{\cal N}_{(\bar{A}, \bar{B}, \bar{C})}}}
\nc{\nabcbt}{{\tilde{\cal N}_{(\bar{A}, \bar{B}, \bar{C})}}}

\nc{\bl}{{\bar{\lambda}}}
\nc{\R}{{\bf R}}
\nc{\C}{{\bf C}}
\nc{\Z}{{\bf Z}}
\nc{\N}{{\bf N}}
\nc{\linf}{{l_{\infty}}}
\nc{\Fb}{{{F} }}
\nc{\tb}{{\bar{t}}}
\nc{\tub}{{\underline{t}}}
\nc{\umguy}{{\| u_m \|_{\infty}}}

\nc{\stable}   {\mbox{SJVL-stable}}
\nc{\Rica}     {{R\overline{ic}}}
\nc{\Ricb}     {\overline{R\overline{ic}}}
\nc{\B}        {{\cal B} }
\nc{\Ltplusd}  {{\cal L}_2  [0, \infty ) }
\nc{\Ltplus}   {{\cal L}_{2+} }
\nc{\Pplus}    {P_{+} }
\nc{\Ltminusd} {{\cal L}_2 (-\infty ,0] }
\nc{\Ltminus}  {{\cal L}_{2-} }
\nc{\Pminus}   {P_{-} }
\nc{\Hf}       {{\cal H}_{\infty} }
\nc{\Ht}       {{\cal H}_{2} }
\nc{\Htp}      {{\cal H}_{2}^{\perp} }
\nc{\RHt}      {{\cal RH}_{2} }
\nc{\RHtp}     {{\cal RH}_{2}^{\perp} }
\nc{\RHf}      {{\cal RH}_{\infty} }
\nc{\Lf}       {{\cal L}_{\infty} }
\nc{\Lt}       {{\cal L}_2 }
\nc{\RLf}      {{\cal RL}_{\infty} }
\nc{\Xminus}   {{\cal X}_- }
\nc{\Xplus}    {{\cal X}_+ }
\nc{\be}       {\begin{equation}}
\nc{\ee}       {\end{equation}}
\nc{\bea}      {\begin{eqnarray}}
\nc{\eea}      {\end{eqnarray}}
\nc{\ba}       {\begin{array}}
\nc{\ea}       {\end{array}}
\nc{\p}        { ^\prime }
\nc{\Xf}       { X_{\infty} }
\nc{\Yf}       { Y_{\infty} }
\nc{\Ymf}      { Y_{tmp } }
\nc{\Hamf}     { H_{\infty} }
\nc{\Sf}       { S_{\infty}}         % PAI 061089
\nc{\Tf}       { T_{\infty}}         % PAI 061089
\nc{\Ff}       { F }
\nc{\Jf}       { J_{\infty} }
\nc{\Jmf}      { J_{tmp } }
\nc{\Ax}       { A_{X} }
\nc{\AFf}      { A_{F} }
\nc{\AFt}      { A_{F_0} }
\nc{\ALt}      { A_{L_2} }
\nc{\CoFf}     { C_{1F} }
\nc{\CoFt}     { C_{1F_0} }
\nc{\BoLt}     { B_{1L_2} }
\nc{\CFt}      { C_{F_2} }
\nc{\sU}       { {\cal U} }
\nc{\Dp}       {D_{\perp}}
\nc{\tR}       {\tilde{R}}
\nc{\tD}       {\tilde{D}}
\nc{\tDp}      {\tD_{\perp}}
\nc{\Doo}      {D_{11}}
\nc{\Dto}      {D_{21}}
\nc{\Dbo}      {D_{\bullet1}}
\nc{\Dob}      {D_{1\bullet}}
\nc{\Fl}       {{\cal F}_{\ell}}
\nc{\fFl}      {\Fl(G,K)}
\nc{\fFlf}     {\|\fFl\|_{\infty}}
\nc{\Gtmp}     {G_{\rm tmp}}
\nc{\Jtmp}     {J_{\rm tmp}}
\nc{\Xtmp}     {X_{\rm tmp}}
\nc{\Ytmp}     {Y_{\rm tmp}}
\nc{\Ftmp}     {F_{\rm tmp}}
\nc{\Ffo}      {F_1}
\nc{\Fft}      {F_2}
\nc{\Fo}       {F^o}
\nc{\Lo}       {L^o}
\nc{\Tot}      {(T_1^{\prime}T_1)^{-1}}
\nc{\Bottil}   {(B_{12}-B_{22}D_{12}^{\prime}D_{11})}
\nc{\Bootil}   {(B_{12}-B_{22}D_{12}^{\prime}D_{11})}
\nc{\Bottilp}  {(B_{12}^{\prime}-\Doo^{\prime}\Dot B_{22}^\prime)}
\nc{\Bootilp}  {(B_{11}^{\prime}-\Doo^{\prime}\Dot B_{21}^\prime)}
\nc{\DoDpDp}   {\Doo^{\prime}\Dp\Dp^{\prime}}
\nc{\xh}       {\hat{x}}
\nc{\xhd}      {\dot{\hat{x}}}
\nc{\tB}       {\tilde{B}}
\nc{\eqdent}   {\hspace{6ex}}
\nc{\tallarray}{\renc{\arraystretch}{1.15}}
\nc{\normarray}{\renc{\arraystretch}{1.0}}

\nc{\tbth}[6]{
  \left[ \begin{array}{ccc}
       #1 & #2 & #3 \\ #4 & #5 & #6 
       \end{array} \right] }
\nc{\thbt}[6]{
  \left[ \begin{array}{cc}
       #1 & #2  \\ #3 & #4 \\ #5 & #6 
       \end{array} \right] }
\nc{\thbth}[9]{
  \left[ \begin{array}{ccc}
       #1 & #2 & #3 \\ #4 & #5 & #6 \\ #7 & #8 & #9
       \end{array} \right] }
\nc{\tbt}[4]{
  \left[ \begin{array}{cc}
       #1 & #2 \\ #3 & #4
       \end{array} \right] }
\nc{\tbo}[2]{
  \left[ \begin{array}{c}
       #1 \\ #2
       \end{array} \right] }
\nc{\tbf}[8]{
  \left[ \begin{array}{cccc}
       #1 & #2 & #3 & #4 \\
       #5 & #6 & #7 & #8
       \end{array} \right] }
\nc{\thbo}[3]{
  \left[ \begin{array}{c}
       #1 \\ #2 \\ #3
       \end{array} \right] }
\nc{\obf}[4]{
  \left[ \begin{array}{cccc}
       #1 &  #2 &  #3 &  #4
       \end{array} \right] }
\nc{\obfi}[5]{
  \left[ \begin{array}{ccccc}
       #1 &  #2 &  #3 &  #4 & #5
       \end{array} \right] }
\nc{\obs}[6]{
  \left[ \begin{array}{cccccc}
       #1 &  #2 &  #3 &  #4 & #5 & #6
       \end{array} \right] }
\nc{\sbo}[6]{
  \left[ \begin{array}{c}
       #1 \\  #2 \\  #3 \\  #4 \\ #5 \\ #6
       \end{array} \right] }
\nc{\fbo}[4]{
  \left[ \begin{array}{c}
       #1 \\  #2 \\  #3 \\  #4
       \end{array} \right] }
\nc{\fibo}[5]{
  \left[ \begin{array}{c}
       #1 \\  #2 \\  #3 \\  #4 \\ #5
       \end{array} \right] }
\nc{\obt}[2]{
  \left[ \begin{array}{cc}
       #1 & #2
       \end{array} \right] }
\nc{\obth}[3]{
  \left[ \begin{array}{ccc}
       #1 & #2 & #3
       \end{array} \right] }
\nc{\tfmat}[4]{
  \left[ \begin{array}{c|c}
       #1 & #2 \\ \hline
       #3 & #4
       \end{array} \right] }
\nc{\LFT}[6]{
  \begin{picture}(60,42)
    \put(20,25){\framebox(20,16){#1}} % P
    \put(24,5){\framebox(12,10){#2}}  % K
    \put(10,40){\makebox(0,0){#3}}  % z
    \put(6,20){\makebox(0,0){#4}}   % y
    \put(50,40){\makebox(0,0){#5}}  % w
    \put(54,20){\makebox(0,0){#6}}  % u
    \put(20,37){\vector(-1,0){20}}
    \put(60,37){\vector(-1,0){20}}
    \put(50,29){\vector(-1,0){10}}
    \put(10,29){\line(1,0){10}}
    \put(36,10){\line(1,0){14}}
    \put(10,10){\vector(1,0){14}}
    \put(10,10){\line(0,1){19}}
    \put(50,10){\line(0,1){19}}
  \end{picture}}
\nc{\lft}[6]{
  \begin{picture}(30,21)
    \put(10,11){\framebox(10,10){#1}}   % P
    \put(11,0){\framebox(8,6){#2}}              % K
    \put(5,21){\makebox(0,0){#3}}               % z
    \put(5,8){\makebox(0,0)[r]{#4}}             % y
    \put(25,21){\makebox(0,0){#5}}              % w
    \put(25,8){\makebox(0,0)[l]{#6}}    % u
    \put(10,19){\vector(-1,0){10}}
    \put(30,19){\vector(-1,0){10}}
    \put(24,13){\vector(-1,0){4}}
    \put(6,13){\line(1,0){4}}
    \put(19,3){\line(1,0){5}}
    \put(6,3){\vector(1,0){5}}
    \put(6,3){\line(0,1){10}}
    \put(24,3){\line(0,1){10}}
  \end{picture}}
\nc{\LFTcenter}[6]{
  \begin{center}
    \LFT{#1}{#2}{#3}{#4}{#5}{#6}
  \end{center}}
\nc{\lftcenter}[6]{
  \begin{center}
    \lft{#1}{#2}{#3}{#4}{#5}{#6}
  \end{center}}
\nc{\lfte}[7]{
  \begin{picture}(80,21)
    \put(0,0){\lft{#1}{#2}{#3}{#4}{#5}{#6}}
    \put(50,10){#7}
  \end{picture}}
\nc{\lftecrunch}[7]{
  \begin{picture}(80,21)
    \put(0,0){\lft{#1}{#2}{#3}{#4}{#5}{#6}}
    \put(40,10){#7}
  \end{picture}}
\nc{\lftss}[6]{
  \begin{picture}(64,42)
    \put(14,25){\framebox(32,16){#1}}
    \put(24,5){\framebox(12,10){#2}}
    \put(8,40){\makebox(0,0){#3}}
    \put(6,20){\makebox(0,0){#4}}
    \put(52,40){\makebox(0,0){#5}}
    \put(54,20){\makebox(0,0){#6}}
    \put(14,37){\vector(-1,0){14}}
    \put(60,37){\vector(-1,0){14}}
    \put(50,29){\vector(-1,0){4}}
    \put(10,29){\line(1,0){4}}
    \put(36,10){\line(1,0){14}}
    \put(10,10){\vector(1,0){14}}
    \put(10,10){\line(0,1){19}}
    \put(50,10){\line(0,1){19}}
  \end{picture}}
\nc{\lftsscenter}[6]{
  \begin{center}
    \lftss{#1}{#2}{#3}{#4}{#5}{#6}
  \end{center}}

\newtheorem{theorem}    {Theorem}
\newtheorem{prop}{Proposition}

\newtheorem{remark}     {Remark}

\begin{document}

\title{\Large \bf {Classical $d$-Step-Ahead Adaptive Control Revisited:
Linear-Like Convolution Bounds and Exponential Stability (Extended Version)}
\author{
Daniel E. Miller and Mohamad T. Shahab\footnote{This research was supported
by a grant from the Natural Sciences Research Council of Canada.} \\
Dept.\ of Elect.\ and Comp.\ Eng., University of Waterloo  \\
Waterloo, ON,
Canada \ N2L 3G1 \\
(miller@uwaterloo.ca) 
}}

%\thispagestyle{empty}
%\date{\today}
%\maketitle
%\thispagestyle{empty}

\maketitle
\thispagestyle{empty}
\pagestyle{empty}

\begin{abstract}
\footnote{This research was supported by the Natural Sciences and Engineering
Research Council of Canada.}
Classical discrete-time adaptive controllers provide asymptotic stabilization
and tracking;
neither exponential stabilization nor a bounded noise gain is typically proven.
In recent work it has been shown, in both the pole placement stability
setting and the first-order one-step-ahead tracking setting,
that if the original, ideal, Projection Algorithm is used
(subject to the common assumption that the plant parameters lie in a convex, compact set and that the parameter estimates are restricted to that set)
as part of the adaptive controller, then
a linear-like convolution bound on the closed loop behaviour can be proven;
this immediately confers exponential stability and a bounded noise gain, 
and it can be leveraged to 
provide
tolerance to unmodelled dynamics and plant parameter variation.
In this paper we extend the approach to the 
$d-$step-ahead adaptive controller 
setting and prove comparable properties.
\end{abstract}

%\noindent
%{\bf Keywords:}
%Adaptive control, Projection algorithm, Exponential stability, Bounded gain.

%\vspace{0.3cm}

\section{introduction}

Adaptive control is an approach used to deal with systems with uncertain
and/or time-varying parameters.
In the classical approach to adaptive control, one combines a
linear time-invariant (LTI)
compensator together with a tuning mechanism to adjust the compensator
parameters to match the plant.
The first general proofs that parameter adaptive controllers could work
came around 1980, e.g. see \cite{morse1978},
\cite{morse1980}, \cite{goodwin1980}, \cite{Narendra1980_pt2}, and
\cite{Narendra1980}. However,
such controllers are typically not robust to unmodelled dynamics, do
not tolerate
time-variations well, have poor transient behaviour,
and do not handle noise (or disturbances) well, e.g. see \cite{rohrs}.
During the following two decades a good deal of research was
carried out to
address these shortcomings. The most common approach was to make
small controller design changes,
such as the use of
signal normalization, deadzones, and $\sigma-$modification,
e.g. see \cite{rick}, \cite{rick2},
\cite{Tsakalis4}, \cite{kreiss}, \cite{Ioa86}.
It turns out that simply using projection
(onto a convex set of admissible parameters)
has proved quite powerful,
and
the resulting controllers typically provide a bounded-noise bounded-state
property, as well as tolerance of some degree of unmodelled dynamics and/or time-variations,
e.g. see \cite{ydstie}, \cite{ydstie2}, \cite{naik}, \cite{wenhill}, \cite{wen} and \cite{hanfu}.
However, in general these controllers provide only asymptotic
stability and not 
exponential stability, with no bounded gain on the noise.
Our goal is to investigate the redesign of adaptive controllers
so that they have more desireable properties.

%The two classical parameter estimators are the Projection Algorithm
%and the Least Squares Algorithm; only the former has the facility
%to handle plant parameter variation.
Here we return to a common approach in classical adaptive control - 
the use of a Projection
Algorithm based estimator together with a tuneable
compensator whose parameters are chosen via the Certainty Equivalence Principle.
In the literature it is the norm to use a modified version of the ideal Projection
Algorithm in order to avoid division by zero;
\footnote{An exception is the work of Ydstie \cite{ydstie}, \cite{ydstie2},
who considers the ideal Projection Algorithm as a special case; however,
a crisp bound on the effect
of the initial condition and a convolution bound on the effect
of the exogenous inputs are not proven. Another notable exception is the 
work of Akhtar and Bernstein \cite{akhtar}, where they are able
to prove Lyapunov stability; however, they do not prove a
convolution bound on the effect
of the exogenous inputs either, and they assume that the high frequency
gain is known.}
it turns out that
an unexpected consequence of this minor adjustment is that some inherent
properties of the scheme are destroyed.
In earlier work by the first co-author
on the
first order setting \cite{scl} and in the pole placement setting of \cite{ccta} and
\cite{rw},
{\bf linear-like 
convolution bounds
on the closed-loop behaviour are proven; such bounds
are highly desirable and have
never before been proven
in the adaptive setting.}
They
confer
exponential stability and a bounded gain on the noise, and
allows a modular approach to analyse
robustness and tolerance to time-varying
parameters.
The objective of the present paper is to use this approach to analyse
the $d-$step-ahead adaptive control problem. While we initially
expected it to follow in a straight-forward manner
from the pole placement setting of \cite{ccta} and \cite{rw},
this has not proven to be the case;
the difficulty stems from the fact that the 
importance of the system delay in this setting 
creates significant additional complexity, as does the fact that in this
problem there is a tracking objective which is not present in 
the pole placement problem.
We have adopted
ideas from \cite{ccta} and \cite{rw} as a starting
point, and we have proven the same
highly desirable linear-like properties
enjoyed in the adaptive pole placement setting.

Before proceeding we present some mathematical preliminaries.
Let $\Z$ denote the set of integers,
$\Z^+$ the set of non-negative integers, $\N$
the set of natural numbers, $\R$ the set of real numbers, and
$\R^+$ the set of
non-negative real numbers.
%We let $\D^0$ denote the open unit disk of the complex plane.
We use the Euclidean $2$-norm for vectors and the corresponding induced
norm for matrices, and denote the norm of a vector or matrix
by $\| \cdot \|$.
We let $\linf ( \R^n)$ denote the set of $\R^n$-valued bounded sequences.
%;
%we define the norm of $u \in \linf ( \R^n)$
%by $\| u \|_{\infty} := \sup_{k \in \Z} \| u(k) \|$.

If ${\cal S} \subset \R^p$ is a convex and compact set, we
define $\| {\cal S} \| := \max_{x \in {\cal S} } \| x \|$ and the
function $\pi_{\cal S} : \R^p \rightarrow
{\cal S}$ denotes the projection onto ${\cal S}$; it is well-known that
$\pi_{\cal S}$
is well-defined.
%With $\eps > 0$ and $c>0$, we let $s ( {\cal S} , c_0 , \eps )$ denote the set of sequences in
%$x \in \linf ( \R^n) $
%taking values in ${\cal S}$ and
%satisfying, for every $t_2 \in \Z$:
%\[
%\sum_{t=t_1}^{t_2-1} \| x(t+1) - x (t) \| \leq c_0 + \eps ( t_2 - t_1 ) , \;
%t_1 > t_2 .
%\]

\section{The Setup}

In this paper
we start with a linear time-invariant discrete-time plant
described by
\begin{eqnarray}
\sum_{i=0}^{n} a_{i} y(t-i) =
\sum_{i=0}^{m} b_{i} u(t-d-i) + w(t),\; t \in \Z,
\label{plant}
\end{eqnarray}
with
\begin{itemize}
\item $y(t) \in \R$ the measured output,
\item
$u(t) \in \R$ the control input,
\item
$w(t) \in \R$ the disturbance (or noise) input;
\item
the parameters regularized so that
$a_0 = 1$,
and 
\item
the system delay is exactly $d$, i.e.
$b_0 \neq 0 $.
\end{itemize}
Associated with this plant model are the polynomials
$A(z^{-1} ) := \sum_{i=0}^n a_i z^{-i}$ and
$ B(z^{-1}) := \sum_{i=0}^{m} b_i z^{-i} $, 
%1 + a_1 z^{-1} + \cdots + a_n z^{-n}  ,
%\]
%\[
%B(z^{-1} ) :=  b_0 + b_1 z^{-1} + \cdots + b_m z^{-m} ,
%\]
as well as the transfer function $z^{-d} \frac{B(z^{-1} )}{A(z^{-1} )}$ and
the list of plant parameters:
\[
\theta_{ab}^* := \left[ \begin{array}{cccccc} {a_1} & \cdots & a_n& 
b_0 & \cdots & b_m \end{array} \right]^T.
\]
It is assumed that $\theta_{ab}^*$ lies in a known set
${\cal S}_{ab} \subset \R^{n+m+1}$.

\begin{remark}
It is straight-forward to verify that if the system has a disturbance at both the input and output, then it can be converted to a system of the above form. 
\end{remark}

The goal is closed-loop stability and asymptotic tracking of
an {\bf exogenous reference input} $y^*(t)$. 
We impose several assumptions on the set of admissible parameters.

\noindent
\framebox[85mm][l]{
\parbox{82mm}{
{\bf Assumption 1:}
The parameter set 
${\cal S}_{ab} $ is compact, and for each $\theta \in {\cal S}_{ab}$, 
the corresponding polynomial
$B(z^{-1} )$
\begin{itemize}
\item
has all of its zeros in the open unit disk, and
\item
the sign of $b_0$ is always the same.
\end{itemize}
}}

\begin{remark}
We have implicitly assumed knowledge of the system delay $d$ as well
as upper bounds on the order of $A(z^{-1})$ and $B(z^{-1})$.
\end{remark}

The boundedness requirement on ${\cal S}_{ab} $ is quite
reasonable in practical situations; it is used here to 
prove uniform bounds and decay rates on the closed-loop behaviour.
The constraint on the zeros of $B(z^{-1} )$ is a requirement that the plant
be minimum phase; this is necessary to ensure tracking of an arbitrary
bounded reference signal \cite{ld_paper}.
Knowledge of the sign of $b_0$ is a common one in 
adaptive control \cite{goodwinsin}.

%The main goal here is to prove a form of stability, with a secondary
%goal that of asymptotic tracking
%of the 
%exogenous reference signal $y^* (t)$.
To proceed we use a parameter estimator together with an adaptive
$d-$step-ahead 
control law.
To design the estimator it is convenient to put the plant into
the so-called predictor form. 
To this end, following \cite{goodwinsin}, we carry out long division by
dividing $A(z^{-1} )$ into one, and define 
$F(z^{-1}) = \sum_{i=0}^{d-1} f_i z^{-i} $ and
$ G( z^{-1} ) = \sum_{i=0}^{n-1} g_i z^{-i} $
satisfying
\[
\frac{1}{A (z^{-1})} = F(z^{-1}) + z^{-d} \frac{ G(z^{-1})}{A (z^{-1} ) } .
\]
Hence, if we define
\[
\beta (z^{-1})  = \sum_{i=0}^{m+d-1} \beta_i z^{-i} :=  F(z^{-1}) 
B( z^{-1} ) ,
\]
\[
\alpha (z^{-1} ) =  \sum_{i=0}^{n-1} \alpha_i z^{-i} :=  G(z^{-1}) ,
\]
\[
\bar{w} (t) := f_0 w (t+d )
 + \cdots + f_{d-1} w ( t+1 ),
\]
then we can rewrite the plant model as
\begin{eqnarray}
y(t+d) &=& \sum_{i=0}^{n-1} \alpha_i y(t-i) + \sum_{i=0}^{m+d-1} \beta_i u(t-i) +
\bar{w} (t) \nonumber 
\end{eqnarray}
\be
=
\underbrace{
\left[ \begin{array}{c}
y(t) \\
\vdots \\
y(t-n+1) \\
u(t) \\
\vdots \\
u( t - m -d +1 ) 
\end{array} \right]^T}_{=: \phi (t)^T}
\underbrace{
\left[ \begin{array}{c}
 \alpha_0 \\
\vdots \\
 \alpha_{n-1} \\
\beta_0 \\
\vdots \\
\beta_{ m +d -1 } 
\end{array} \right]^T}_{=: \theta^*} 
 + \bar{w} (t), \;\; t \in \Z . 
\label{newplant}
\ee

Let ${\cal S}_{\alpha  \beta}$ denote the set of 
admissible $\theta^*$ which arise from the original plant parameters
which lie in ${\cal S}_{ab}$;
since the associated mapping is continuous, it is clear that
the compactness of 
${\cal S}_{ab}$ 
means that
${\cal S}_{\alpha  \beta}$ 
is compact as well.
Furthermore, it is easy to see that $f_0 =1$, so $\beta_0 = b_0$, which
means that
the sign of $\beta_0$ is always the same.
It is convenient that the set of admissible parameters in the
new parameter space be convex and closed; so at this point let
${\cal S} \subset \R^{n+m+d}$ be any {\bf compact and convex set} containing
${\cal S}_{\alpha  \beta}$ for which the $n+1^{th}$ element  (the one
which corresponds to $\beta_0$)
is never zero, e.g. the convex hull of ${\cal S}_{\alpha  \beta}$ would do.

The $d-$step-ahead control law is the one given by
\[
y^*(t+d) = 
\sum_{i=0}^{n-1} \alpha_i y(t-i) + \sum_{i=0}^{m+d-1} \beta_i u(t-i) ;
\]
in the absence of a disturbance, and assuming that this controller
is applied for all $t \in \Z$,
we have
$y (t) = y^* (t)$ for all $ t \in  \Z$.
Of course, if the plant parameters are unknown, we need to
use estimates; also, the adaptive version of the
$d$-step-ahead control law is only applied after some
initial time, i.e. for  $ t \geq t_0$.

\subsection{Initialization}

In most adaptive controllers the goal is to prove asymptotic results,
so the details of the initial condition is unimportant.
Here, however, we wish to get a bound on the transient behaviour so we must
proceed carefully.
%In the pole placement setting of \cite{ccta}, this was relatively
%straight-forward:
%the delay plays no role,
%the controller is strictly
%causal, and
%we start the 
%plant estimator off at time $t_0$, with an initial "plant state" of 
%$\phi (t_0) = \phi_0$.
%Here we have a more complicated situation, even in
%the case of $d=1$, since the
%proposed controller is not strictly causal.
%
%To proceed,
%first let us examine the predictor form (\ref{newplant}) with
%$t$ replaced by $t-d$: if we apply the $z-$transform, we easily see that
%there are initial conditions of
If we wish to solve (\ref{newplant}) for $y(t)$ starting at
time $t_0$, it is clear that we need an initial condition of
\[
x_0 :=
 \left[ \begin{array}{ccc}
 y( t_0 -1) & \cdots & y(t_0 -n-d+1)  \end{array} \right.
\]
\[
\;\;\; \left. \begin{array}{ccc}  u(t_0 -1 ) & \cdots &
u( t_0 -m-2d+1 ) \end{array} \right]^T .
\]
%this is sufficient information to obtain
%$\{ \phi (t_0-1 ) ,..., \phi (t_0 -d) \}$.

%\begin{remark}
%Observe that if $d \neq 1$, the $d-$step-ahead model is not minimal, 
%whereas the initial model (\ref{plant}) is minimal if $A (z^{-1})$ and
%$B(z^{-1})$ is coprime. 
%Hence, we could get away with fewer initial conditions,
%\end{remark}

\subsection{Parameter Estimation}

We can rewrite the plant (\ref{newplant}) as
\be
y(t+1) = \phi ( t-d +1 ) ^T \theta^* +  \bar{w} (t-d+1) , \; t \geq t_0-1.
\label{newplantn}
\ee
Given an estimate $\hat{\theta} (t)$ of $\theta^*$ at time $t$,
we define
the {\bf prediction error} by
\[
e(t+1) := y(t+1) -  \phi (t-d+1)^T  \hat{\theta} (t) ;
\]
this is a measure of the error in $\hat{\theta} (t)$.
A common way to
obtain
a new estimate is
from
the solution of the optimization problem
\[
argmin_{\theta} \{ \| \theta - \hat{\theta} (t) \| : y(t+1) = \phi (t-d+1)^T {\theta} \} ,
\]
yielding the
{\bf ideal (projection) algorithm}
\[
\that (t+1) =
\]
\be
\left\{
\begin{array}{ll}
\that (t) & \mbox{ if $\phi (t-d+1) = 0 $} \\
\that (t) + \frac{\phi (t-d+1)}{\| \phi (t-d+1)\|^2}
\,  e (t+1)  & \mbox{ otherwise;}
\end{array}
\right.
\label{orig}
\ee
at this point, we can also restrain it to ${\cal S}$ by projection.
Of course, if $\phi (t-d+1)$ is close to zero, numerical problems can occur, so it is the norm
in the literature (e.g. \cite{goodwin1980} and \cite{goodwinsin}) to add a constant to
the denominator, but as pointed out in \cite{scl}, \cite{ccta}, and \cite{rw}, this can lead to the loss
of exponential stability and a loss of a bounded gain on the noise. We propose a middle
ground: as proposed in \cite{ccta} and \cite{rw}, we turn off the 
estimation if it is clear than the disturbance signal
$\bar{w} (t)$ is swamping the estimation error.
To this end, with
$\delta \in (0, \infty ]$, we turn off the estimator if the
update is larger than $2 \| {\cal S} \| + \delta$ in magnitude; so define
$\rhod : \R^{n+m+d} \times \R  \rightarrow \{ 0,1 \}$ by
\[
\rhod ( \phi (t-d+1) , e (t+1)  ) :=
\]
\[
\left\{ \begin{array}{ll}
1 & \;\; \mbox{ if } | e(t+1) | < ( 2  \| {\cal S} \| + \delta )
\| \phi (t-d+1) \| \\
0 & \;\; \mbox{otherwise; }
\end{array} \right.
\]
given $\hat{\theta} (t_0 -1) = \theta_0$, for $t \geq t_0-1$ we define\footnote{If $\delta = \infty$, then we adopt the understanding that $\infty \times 0 = 0$, in which case this formula collapses into the original one.}
\[
\check{\theta} (t+1) = \hat{\theta} (t) + \rhot \times
\]
\be
\;\;\;\;
\frac{ \phi (t-d+1)} {\| \phi (t-d+1) \|^2} e(t+1) ,
\label{estimator1}
\ee
which we then
project onto ${\cal S}$:
\be
\hat{\theta} (t+1):= \pi_{\cal S} ( \check{\theta} (t+1) ).
\label{estimator2}
\ee

\subsection{Properties of the Estimation Algorithm}

Analysing the closed-loop system will require a careful
analysis of the estimation algorithm.
We define
the
parameter estimation error by
$\tilde{\theta} (t): = \hat{\theta} (t) - \theta^* $
and the corresponding Lyapunov function associated with
$\tilde{\theta} (t)$, namely
$V(t) : = \tilde{\theta} (t)^T \tilde{\theta} (t) $.
In the following result we list a property of $V(t)$; it is a straight-forward
generalization
of what holds in the pole placement setup of \cite{ccta} and \cite{rw}.
%
%\noindent
%\framebox[160mm][l]{
%\parbox{155mm}{
\noindent
\framebox[85mm][l]{
\parbox{82mm}{
\begin{prop}
For every $t_0 \in \Z$, $\phi_0 \in \R^{n+m+d}$, ${\theta}_0 \in {\cal S}$,
$\theta_{ab}^* \in {\cal S}_{ab}$,
$y^*, w \in \linf$, and $\delta \in ( 0, \infty ]$, when the estimator 
(\ref{estimator1}) and (\ref{estimator2})
is applied to the plant (\ref{plant}),
the following holds:
\[
 \| \hat{\theta} (t+1) - \hat{\theta}  (t) \|
\leq   \rhot \times
\]
\be
\;\; \frac{ |e(t+1)| }{ \|  \phi (t-d+1) \| } ,
\; t \geq t_0-1 ,
\label{errorsum}
\ee
\[
V(t)
 \leq 
V( t_0-1 ) +   \sum_{j=t_0-1}^{t-1}
\rhoj \times
\]
\[
\;\;\;
[ -\frac{1}{2}
\frac{[e (j+1) ]^2}{ \| \phi (j-d+1) \|^2} +
2 \frac{[ \bar{w}(j-d+1)]^2}{ \| \phi (j-d+1) \|^2}] , \;\; t \geq t_0 -1.
\]
\end{prop}
}}

\subsection{The Control Law}
The elements of $\hat{\theta} (t)$ are partitioned in a natural way as
%\[
%\hat{\theta} (t) :=
%\]
\[
\left[
\begin{array}{cccccc}
\hat{\alpha}_0 (t)  &
\cdots &
\hat{\alpha}_{n-1} (t) &
 \hat{\beta}_0 (t) &
\cdots &
 \hat{\beta}_{m+d-1} (t) 
\end{array}
\right]^T .
\]
The {\bf one-step-ahead adaptive control law} is that of
\[
y^*(t+d) = \hat{\theta} (t)^T \phi (t) ,
\; t \geq t_0 ,
\]
or equivalently
\be
\sum_{i=0}^{m+d-1} \hat{\beta}_{i} (t) u(t-i) = 
y^*(t+d) -
 \sum_{i=0}^{n-1} \hat{\alpha}_{i} (t) y(t-i) .
\label{controller}
\ee
Hence, as is common in this setup, 
we assume that the controller has access to the reference signal
$y^*(t)$ exactly $d$ time units in advance.

\begin{remark}
With this choice of control law, it is easy to prove that the prediction
error $e(t)$ and the {\bf tracking error}
\[
\eps (t) := 
y^* (t) - y(t) 
\]
are different if $d \neq 1$. Indeed, 
it is easy to see that
\be
\eps (t) = -\phi (t-d)^T \tilde{\theta} (t-d) + \bar{w} (t-d) ,
\; t \geq t_0 + d, 
\label{epsdude}
\ee
\be
e(t) = - \phi (t-d)^T \tilde{\theta} (t-1) + \bar{w} (t-d) ,
\; t \geq t_0.
\label{edude}
\ee
%\be
%\eps(t) 
%&=& 
%e(t) + \phi(t-d)^T [ \tilde{\theta} (t-1) - \tilde{\theta} (t-d) ] \\
%=
%e(t) + \phi(t-d)^T [ \hat{\theta} (t-1) - \hat{\theta} (t-d) ] ,
%\; t \geq t_0.
%\label{eedude}
%\ee
\end{remark}

%At this point construct two closed-loop models. The first one will be useful
%for analysing the closed-loop behaviour.
{\bf The goal of this paper is to prove that the adaptive controller} consisting of
the estimator 
(\ref{estimator1})-(\ref{estimator2}) together with the control equation
(\ref{controller}) {\bf yields highly desirable linear-like convolution
bounds on the closed-loop behaviour}.
While the approach is similar to that in our earlier work \cite{ccta} and \cite{rw},
it requires a much more nuanced analysis.
In the next section we develop several models used in the development, after which we
state and prove the main result.

\section{Preliminary Analysis}

\subsection{A Good Closed-Loop Model}

In our pole-placement adaptive control setup \cite{ccta}, \cite{rw}, a key 
closed-loop model consists of an update equation for $\phi (t)$, with the
state matrix consisting of controller and plant estimates;
this was effective - the characteristic polynomial of this matrix
is time-invariant and has all roots in the open unit disk.
If we were to mimic this in the
one-step-ahead setup, the characteristic polynomial would have
roots which are time-varying, with some at zero and the rest at the zeros of
$\hat{\beta} (t, z^{-1})$, which is time-varying and may not have roots in the open unit disk.
Hence, at this point we make an important deviation from the approach
of \cite{ccta} and \cite{rw} and construct the following
update equation for $\phi (t)$ which avoids the use of plant parameter estimates,
but is driven by the tracking error.
Only two elements of $\phi$ have a complicated
description:
\begin{eqnarray*}
\phi_1 (t+1) &=& y(t+1) 
= \eps(t+1) + y^*(t+1) , \\
&& \;\;\; t \geq t_0+d-1,
%& = & y^*(t+1) +
%e(t+1) + \\
%&& \underbrace{\phi(t-d+1)^T [ \tilde{\theta} (t) - \tilde{\theta} (t-d+1) ]}_{=: \delta_1 (t)} 
\end{eqnarray*}
and
the $u(t+1)$ term, for which we use the original plant model to write:
\[
\phi_{n+1}(t+1) = u(t+1) 
\]
\begin{eqnarray*}
%&=&
% \frac{1}{b_0} [ \sum_{i=0}^d a_i y(t+d+1-i) + \\
%&&
%\sum_{i=d+1}^n a_i y(t+d+1-i) - \\
%&&   \sum_{i=0}^{m-1} b_{i+1} u(t-i) - w ( t+d+1 ) ]\\
&=& 
\frac{1}{b_0} [ \sum_{i=0}^d a_i ( \eps (t+d+1-i) + y^*  (t+d+1-i) ) + 
\end{eqnarray*}
\[
 \sum_{i=d+1}^n a_i y(t+d+1-i) - \\
   \sum_{i=0}^{m-1} b_{i+1} u(t-i) - w ( t+d+1 ) ] .
\]
%\[
%\underbrace{ \phi ( t+1)^T [ \hat{\theta} ( t+d ) - \hat{\theta} (t+1) ] +
%\alpha_0 \phi (t)^T [  \hat{\theta} ( t+d-1 ) - \hat{\theta} (t) ] }_{=:
%\delta_2 (t)} .
%\]
With $e_i \in \R^{n+m+d}$ the $i^{th}$ normal vector,
if we now define
\be
B_1:= 
e_1 , 
%\left[ \begin{array}{c}
%1 \\ 0 \\ \vdots \\ 0 \end{array} \right] ,
\;\;
B_2 := e_{n+1} ,
%  \left[ \begin{array}{c}
%0 \\ \vdots \\ 0 \\ 1 \\ 0 \\ \vdots \\ 0
%\end{array} \right] ,
\label{rdef}
\ee
then 
it is easy to see that there exists a matrix $A_g \in \R^{(n+m+d) \times
(n+m+d )}$ so that
the following equation holds:
\[
\phi (t+1) = {A}_g  \phi (t) + B_1  \eps (t+1) + 
\]
\begin{eqnarray}
&& 
  B_2 
 \sum_{j=0}^d [ \frac{a_{d-j}}{b_0}  \eps (t+1+j) + 
\frac{a_{d-j}}{b_0}  y^*  (t+1+j) ] \nonumber + \\
&&  B_1 y^*(t+1) -
\frac{1}{b_0} B_2
 w ( t+d+1 ) , \; t \geq t_0-1 . \label{good_model}
\end{eqnarray}
The characteristic equation of $A_g$ 
equals
$\frac{1}{b_0} z^{n+m+d} B( z^{-1} )$, so all of its roots are in the
open unit disk.

\subsection{A Crude Closed-Loop Model}

At times we will need to use a crude model to bound the size of the
growth of $\phi (t)$ in terms of the exogenous inputs.
%Here we build such a model.
Once again, only two elements of $\phi (t)$ have a complicated
description:
to describe $y(t+1)$ we use the plant model:
\begin{eqnarray*}
\phi_1 (t+1) &=& y(t+1)  \nonumber \\
&=& 
- \sum_{i=1}^{n} a_{i} y(t+1-i) +  \nonumber \\
&&  \sum_{i=0}^{m} b_{i} u(t+1-d-i) + w(t+1) \\
&=: & \bar{\theta}_{ab}^* \phi (t) + w(t+1), 
\nonumber
\end{eqnarray*}
and to describe $u(t+1)$ we use the control law:
\begin{eqnarray*}
y^*(t+d) &=& \hat{\theta} (t) ^T \phi (t) \\
\Rightarrow \; y^* (t+d+1) &=& \hat{\theta} (t+1) ^T \phi (t+1) , \; t \geq t_0
-1 ;
\end{eqnarray*}
it is easy to define
$\bar{\theta}_{\alpha \beta} (t)$ in terms of the elements of
$\hat{\theta} (t+1)$ so that
\[
y^* (t+d+1) = \bar{\theta}_{\alpha \beta} (t)^T \phi (t) +
\]
\[
\;\;\; 
\hat{\alpha}_0 (t+1) y(t+1) + \hat{\beta}_0 (t+1) u(t+1) , \; t \geq t_0-1 .
\]
If we combine this with the formula for $y(t+1)$ above, we end up with
\[
u(t+1) = \frac{1}{ \hat{\beta} _0 (t+1)} [
- \bar{\theta}_{\alpha \beta} (t) - \hat{\alpha}_0 (t+1) \bar{\theta}_{ab}^* ]
\phi (t) + 
\]
\[
 \frac{1}{\hat{\beta}_0 (t+1)} y^* (t+d+1) -
\frac{\hat{\alpha}_0 (t+1)}{\hat{\beta}_0 (t+1)} w(t+1) , 
\; t \geq t_0-1 .
\]
Hence, we can define matrices
$A_b (t)$, $B_3 (t)$ and $B_4 (t)$ so that
\begin{eqnarray}
\phi (t+1) &=& A_b (t) \phi (t) + B_3 (t) y^* (t+d+1) + \nonumber \\
&&  B_4 (t) w(t+1), \; t \geq t_0-1; \label{crude}
\end{eqnarray}
due to the compactness of ${\cal S}_{ab}$, ${\cal S}_{\alpha \beta}$ and ${\cal S}$, the following
is immediate:

\noindent
\framebox[85mm][l]{
\parbox{82mm}{
\begin{prop}
There exists a constant $c_1$ so that for every
$t_0 \in \Z$, $\phi_0 \in \R^{n+m+d}$, ${\theta}_0 \in {\cal S}$,
$\theta_{ab}^* \in {\cal S}_{ab}$,
$y^*,w \in \linf$, and $\delta \in ( 0, \infty ]$, 
when the adaptive controller
(\ref{estimator1}), (\ref{estimator2})
and
(\ref{controller}) is applied
to the plant (\ref{plant}),
the following holds:
\[
\| A_b (t) \| \leq c_1 , \; \| B_3 (t) \| \leq c_1 , \; \| B_4 (t) \|
\leq c_1 , \; t \geq t_0-1 .
\]
\end{prop}
}}

\subsection{A Better Closed-Loop Model}

The good closed-loop model (\ref{good_model})
is driven by future tracking error signals.
We can now combine this 
with the crude closed-loop model (\ref{crude}) to create a new model
which is driven by perturbed versions of the present
and past values of $\phi$, with the weights associated with the parameter
update law. To this end, first define
\[
\nu (t-1) := 
\rho_{\delta} ( \phi (t-d) , e(t))
\times
\frac{ \phi (t-d)} {\| \phi (t-d) \|^2} e(t) , \; t \geq t_0.
\]
The following result plays a pivotal role in the analysis of the
closed-loop system.

\noindent
\framebox[85mm][l]{
\parbox{82mm}{
\begin{prop}
There exists a constant $c_2$ so that
for every $t_0 \in \Z$, $\phi_0 \in \R^{n+m+d}$, ${\theta}_0 \in {\cal S}$,
$\theta_{ab}^* \in {\cal S}_{ab}$,
$y^*, w \in \linf$, and $\delta \in ( 0, \infty ]$, 
when the adaptive controller
(\ref{estimator1}), (\ref{estimator2})
and
(\ref{controller}) 
is applied to the plant (\ref{plant}),
the following holds:
\[
\phi (t+1) = A_g \phi (t) + \sum_{j=0}^{d-1}  \Delta_j (t) \phi (t-j) +   \eta (t) ,
\; t \geq t_0+d-1 ,
\]
with
\[
\| \eta (t) \| \leq  c_2 ( 1 + \nu ( t+2) + \cdots + \nu (t+d+1)) \times
\]
\[
\;\; [ \sum_{j=1}^{d+1} | y^* ( t+j) | + 
\sum_{j=1}^{d+1} ( | w(t+j) | +
 | \bar{ w} ( t+j ) | ) ]
\]
and
\[
\| \Delta_j (t) \| \leq c_2 ( \nu ( t-d+2) +  \cdots + \nu (t+d+1)) ,
\]
\[
\;\;\;\;\;\;\;
\; j=0,...., d-1 .
\]
\end{prop}
}}

\vspace{0.3cm}
\noindent
{\bf Proof:} See the Appendix.
$\Box$

To make the model of Proposition 3 amenable to analysis, we define a 
new extended state variable and associated matrices:
\be
\bar{\phi} (t) := \fbo{\phi (t)}{\phi (t-1)}{\vdots}{\phi (t-d+1)} , \;
%\bar{B}_1 := \left[ \begin{array}{c} {I} \\ {0} \\ {\vdots} \\ {0}
%\end{array} \right],
% \in \R^{d(n+m+d) \times {n+m+d}} ,
%\]
%\be
A_{nom} =
\left[ \begin{array}{cccc}
A_g  &  & & \\
I &  & &  \\
  & \ddots  & & \\
 & & I & 0 
\end{array} \right] ,
\label{anom}
\ee
and
\be
\bar{B}_1 := \left[ \begin{array}{c} {I} \\ {0} \\ {\vdots} \\ {0}
\end{array} \right],
\Delta (t) =
\left[ \begin{array}{cccc}
\Delta_0 (t)  & \Delta_1 (t)  & \cdots & \Delta_{d-1} (t)  \\
0 & \cdots  & \cdots  & 0 \\
\vdots  & \cdots & \cdots & \vdots \\
 0& 0  & 0  & 0
\end{array} \right] ,
\label{ddef}
\ee
which gives rise
to a state-space model which will play a key role
in our analysis:
\be
\bar{\phi} (t+1) =
[ A_{nom} + \Delta (t) ] \bar{\phi} (t) + \bar{B}_1 \eta (t) ,
\; t \geq t_0+d-1 .
\label{starplant}
\ee
Now $A_g$ arises from $\theta_{ab}^* \in {\cal S}_{ab}$, and
lies in a corresponding compact set ${\cal A}$; furthermore, its
eigenvalues are at the zeros of $B(z^{-1})$ which has all of its
roots in the open unit disk, so we can use classical arguments
to prove that 
there exists $\gamma_1$ and $\sigma \in (0,1)$ 
so that 
\[
\| A_{nom} ^i \| \leq \gamma_1 \sigma^i , \; i \geq 0 .
\]
Indeed, 
we can choose any $\sigma$ larger than
\[
\underline{\lambda} :=
\max_{\theta_{ab}^* \in {\cal S}_{ab}} \{| \lambda| : B( \lambda^{-1} ) =0 \} .
\]

Equations of the form given in 
(\ref{starplant}) arise in classical adaptive control approaches; 
the following proposition follows easily from
the lemma of Kreisselmeier \cite{kreiss2}.
%{Furthermore, in
%\cite{kreiss2} it is assumed that $\alpha_i $ and
%$\beta_i $ are strictly greater than zero, but it is trivial to extend this
%to allow for zero as well.}

\noindent
%\framebox[160mm][l]{
%\parbox{155mm}{
\noindent
\framebox[85mm][l]{
\parbox{82mm}{
\begin{prop} 
Consider the discrete-time system
(\ref{starplant})
with $\Phi (t, \tau )$ denoting the 
state transition matrix corresponding to $A_{nom} + \Delta (t)$.
Suppose that there exist constants
$\beta_i \geq 0$
so that 
for all $t > \tau \geq t_0+d-1 $ we have
\[
\sum_{i = \tau}^{t-1} \|  \Delta (i) \| \leq
\beta_0 + \beta_1 ( t- \tau )^{1/2} + \beta_2 ( t- \tau ) ,
\]
and there exists a $\mu \in ( \sigma , 1 )$ and $N \in \N$ satisfying
\[
 {\beta_2} < \frac{1}{ \gamma_1} ( \frac{\mu}{\gamma_1^{1/N}} -
\sigma ) .
\]
Then
there exists a constant $\gamma_2$ so that 
the transition matrix satisfies
\[
\| \Phi (t, \tau  ) \| \leq \gamma_2 \mu^{t - \tau }  , \; t \geq \tau .
\]
\end{prop}
}}

%\vspace{0.2cm}
%\noindent
%{\bf Proof:} 

%In \cite{kreiss2} this is proven, in particular, for the case of a fixed
%$A_{nom}$. Indeed, formulas are 
%provided for $\gamma_2$ and $\mu$ in terms of
%$\gamma_1$, $\sigma$, $\beta_0$, $\beta_1$ and $\beta_2$. It follows 
%immediately that the same results hold for our situation, in which
%$A_{nom}$ lies in a compact set ${\cal A}$.
%$\Box$

\vspace{0.2cm}

\section{The Main Result}

\noindent
\framebox[85mm][l]{
\parbox{82mm}{
\begin{theorem}
For every
$\delta \in ( 0 , \infty ]$ and
$\lambda \in ( \underline{\lambda} ,1)$ there exists a constant $c>0$ so that for
every $t_0 \in \Z$, 
plant parameter ${\theta}_{ab}^*  \in {\cal S}_{ab}$,
% (which maps
%to $\theta^* \in S_{\alpha \beta} \subset S$),
exogenous signals $y^*, w \in \ell_{\infty}$,
estimator
initial condition $\theta_0 \in {\cal S}$,
and plant initial condition
\[
x_0 = \left[ \begin{array}{ccc}
 y( t_0 -1) & \cdots & y(t_0 -n-d+1)  \end{array} \right.
\]
\[
\;\;\; \left. \begin{array}{ccc}  u(t_0 -1 ) & \cdots & 
u( t_0 -m-2d+1 ) \end{array} \right]^T ,
\]
when the adaptive controller 
(\ref{estimator1}), (\ref{estimator2})
and
(\ref{controller}) is applied
to the plant (\ref{plant}), the following bound holds:
\[
\| \phi (k) \| \leq c \lambda^{k- t_0} \| x_0 \| +
\]
\be
\sum_{j=t_0}^{k}  c \lambda^{k-j} ( | y^* (j) | + | w(j)| ) , \;\;
k \geq  t_0 ;
\label{thm1_eq}
\ee
furthermore, if $w =0$ then
\[
\sum_{k=t_0+2d-1}^{\infty} \eps (k)^2 \leq c ( \| x_0 \|^2 + \sup_{j \geq t_0} |
y^* (j) |^2 ).
\]
\end{theorem}
}}

\begin{remark}
Theorem 1 implies that the system has a bounded gain (from 
$y^*$ and $w$ to $y$) in every $p-$norm.
\end{remark}

\begin{remark}
Most $d-$step-ahead adaptive controllers are 
proven to yield a weak form of stability,
such as boundedness (in the presence of a non-zero disturbance) or asymptotic
stability (in the case of a zero disturbance), which means that details surrounding initial conditions
can be ignored.
Here the goal is to prove a stronger, linear-like, convolution
bound as well as exponential stability, so it requires more detailed analysis.
\end{remark}

\begin{remark}
In the absense of noise, most d-step-ahead adaptive controllers simply say
that the tracking error is square summable, e.g. see \cite{goodwinsin}.
Here we prove something much stronger: we
provide an upper bound on the $2-$norm in terms of the size of the
initial condition and the $\infty$-norm of the reference signal.
\end{remark}

\vspace{0.2cm}
\noindent
{\bf Proof:}
This proof is based on a significant modification to our earlier proof 
of the adaptive pole placement controller \cite{ccta}, \cite{rw}.
The proof is more complicated here for two
reason: first of all, here we have to analyse $\bar{\phi} (t)$ rather than
$\phi (t)$; second of all, Proposition 3 provides a bound
on $\| \Delta (t) \|$ which not only depends on 
$\nu (t+d)$ but also many other values of $\nu ( \cdot)$.

Fix $\delta \in (0, \infty ]$ and $\lambda \in ( \underline{\lambda}, 1)$. Let
 $t_0 \in \Z$,
$\theta_{ab}^* \in {\cal S}_{ab}$,
$y^*, w\in \linf$, $\theta_0 \in {\cal S}$, and $x_0 \in \R^{n+m+3d}$
be arbitrary.
Now choose
$\lambda_1 \in
( \underline{\lambda} ,  \lambda )$.
Observe that $x_0$ gives rise to $\phi (t_0 -1)$,..., $\phi (t_0 - d+1)$, 
as well
as $\bar{\phi} (t_0 -1)$, which we label
$\bar{\phi}_0$; it is clear that
$\| \bar{\phi}_0 \| \leq d \| x_0 \| $ and
\[
\| \phi ( t_0-j) \| \leq \| x_0 \| , \; j=1,..., d-1 .
\]
%\[
%\| \bar{\phi}_0 \| \leq d \| x_0 \| .
%\]
To proceed we will analyse (\ref{starplant})
and obtain a bound on $\bar{\phi} (t)$ in terms of $\eta (t)$, 
$\bar{w} (t)$, and $y^* (t)$, which we will then convert to the desired form.
First of all,
we see from Proposition 3 that there exists a constant
$c_2$ so that
\[
\| \eta (t) \| \leq  c_2 ( 1 + \nu ( t+2) + \cdots + \nu (t+d+1)) \times
\]
\be
\;\; \underbrace{[ \sum_{j=1}^{d+1} | y^* ( t+j) | +
\sum_{j=-2d+1}^{d+1} ( | w(t+j) | 
+ | \bar{ w} ( t) | ) ]}_{=: \tilde{w} (t)}
\label{etabd}
\ee
and
\[
\| \Delta (t) \| \leq c_2  ( \nu ( t-d+2) +  \cdots + \nu (t+d+1)) ,
\]
\be
\;\;t \geq t_0 + d - 1 .
\label{dddbd}
\ee
Before proceeding, we provide a useful preliminary result; it
follows
immediately from Proposition 2.

\vspace{0.1cm}
\noindent
{\bf Claim 1:}
There exists a constant $\gamma_3$ so that
\[
\| \phi (t+i) \| \leq \gamma_3 \| \phi (t) \| +
\gamma_3  \sum_{j=1}^{i} [ | y^* ( t+d+j) | + | w(t+j) | ] 
\]
for
$t \geq t_0-1 $ and
$i=1,..., 2d $.
%and
%\[
%\| \bar{\phi} ( t+d) \| \leq \gamma_3 \| \phi (t) \| +
%\gamma_3  \sum_{j=1}^{d} [ | y^* ( t+d+j) | + | w(t+j) | ]  .
%\]

%\vspace{0.3cm}
%\noindent
%{\bf Proof:} This follows from solving (\ref{crude}) iteratively and
%applying Proposition 2.
%\noindent
%$\Box$

In order to apply Proposition 4, we need to compute a bound on a sum;
the following result provides an avenue.

\vspace{0.3cm}
\noindent
{\bf Claim 2:}
There exists a constant $\gamma_4$ so that
for every $t_2 > t_1 \geq t_0+d-1$, 
\[
\sum_{j=t_1}^{t_2-1} \| \Delta (j) \| \leq
\gamma_4 [ \sum_{j=t_1-d+2}^{t_2+d} \nu(j)^2 ]^{1/2} (t_2- t_1 )^{1/2} .
\]

\vspace{0.3cm}
\noindent
{\bf Proof:} 
It follows from (\ref{dddbd}) that
\[
\sum_{j=t_1}^{t_2} \| \Delta (j) \| \leq
2
c_2 d \sum_{j=t_1-d+2}^{t_2+d} |\nu(j) | .
\]
If we apply the Cauchy-Schwarz inequality 
%we have
%\[
%\sum_{j=t1-d+2}^{t_2+d} | \nu (j) | \leq
%[ \sum_{j=t_1-d+2}^{t_2+d} |\nu(j) | ]^{1/2}
%(t_2- t_1 + 2d-2 )^{1/2} .
%\]
%Since
and observe that
$(t_2-t_1)^{1/2} \leq ( 2d)^{1/2} ( t_2 - t_1 + 2d-2)^{1/2} $,
then the result follows.

\noindent
$\Box$
\vspace{0.3cm}

At this point we consider two cases: the easier case in which there is no noise,
and the harder case in which there is noise.

\vspace{0.3cm}
\noindent
{\bf Case 1}: $w(t) = 0$, $t \geq t_0-d$.

Using the definition of $\nu (j)$,
the bound on $V(t)$ given by Proposition 1
simplifies to
\[
V(t) \leq V ( t_0-1  ) - \frac{1}{2} \sum_{j=t_0-1}^{t-1} \nu (j)^2 ,
 \; t \geq t_0  .
\]
Since $V( \cdot ) \geq 0$ and $ V( t_0-1  ) = \| \theta_0 - \theta^* \|^2
\leq 4 \| {\cal S} \|^2$, this means that
\be
\sum_{j=t_0-1}^{t-1} \nu (j)^2 \leq 8 \| {\cal S} \| ^2, \; t \geq t_0 .
\label{nubd}
\ee
From Claim 2 we conclude that
\[
\sum_{j=t_1}^{t_2-1} \| \Delta (j)  \|  \leq
8^{1/2} \gamma_4 \| {\cal S} \| 
( t_2-t_1)^{1/2} 
 , \; t_2 > t_1 \geq t_0+d .
\]
Now we apply Proposition 4: we set
\[
\beta_0 = \beta_2 = 0 ,
\; \beta_1 = 8^{1/2} \gamma_4 \| {\cal S} \| ,
\; \mu = \lambda.
\]
It is now trivial to choose $N \in \N$ so that
$\frac{\lambda}{ \gamma_1^{1/N}} - \lambda_1 > 0 $,
namely
$N = \mbox{int} [ \frac{\ln ( \gamma_1 )}{\ln ( \lambda ) - \ln ( \lambda_1 ) } ] + 1 $,
which means that
${\beta_2} < \frac{1}{ \gamma_1 } (  \frac{\lambda}{ \gamma_1^{1/N}} - \lambda_1 ) $.
From Proposition 4 we see that there exists a constant $\gamma_2$ so that the state
transition matrix $\Phi ( t, \tau )$ corresponding to ${A}_{nom}  + \Delta (t)$ satisfies
\[
\| \Phi (t, \tau ) \| \leq \gamma_2 \lambda^{t - \tau } , \;
t \geq \tau \geq  t_0 +d .
\]
Also, we see from (\ref{etabd}), (\ref{nubd}) and Proposition 3 that
\[
| \eta (t) | \leq c_2 ( 1 + 8^{1/2} \| S \| ) \tilde{w} (t).
\]
If we now apply this to
(\ref{starplant}), we see that there exists a constant $\gamma_5$ so that
\[
\| \bar{\phi} (k) \| \leq \gamma_5 \lambda^{k- t_0} \| \bar{\phi} ( t_0+d) \| +
\]
\[
\;\;\;
\sum_{j=t_0+d}^{k-1}  \gamma_5 \lambda^{k-1-j} | \tilde{w} (j) |  , \;\;
k \geq  t_0+d .
\]
At this point we can use Claim 1 to find a bound on 
$ \| \bar{\phi} ( t_0+d) \|$ in terms of $x_0$, $y^*$ and $w$; if we convert the bounds on
$| \tilde{w} (j) | $ to bounds on
$| w(t)|$ and $|y^* (t)|$, 
then (\ref{thm1_eq}) holds for this case.

\noindent
{\bf Case 2:} $w(t) \neq 0$ for some $t \geq t_0-d$.

This case is much more involved since noise can radically
affect parameter estimation. Indeed, even if the parameter estimate is
quite accurate at a point in time, the introduction of a large noise signal
(large relative to the size of $\phi (t)$) can create a highly inaccurate
parameter estimate.
Following \cite{ccta} and \cite{rw},
we partition the timeline into two parts: one in which the noise is
small versus $\phi$ and one where it is not.
To this end, with
$\xi >0 $ to be chosen shortly,
partition $\{ j \in \Z : j \geq t_0 \}$ into $S_{good}$ and $S_{bad}$,
respectively:
\[
\{ j \geq t_0 : \phi (j-d+1) \neq 0 \mbox{ and }
 \frac{[ \bar{w}(j-d+1)]^2}{ \| \phi (j-d+1) \|^2}  < \xi \} ,
\]
\[
\{ j \geq t_0 : \phi (j-d+1) = 0 \mbox{ or }
 \frac{[ \bar{w}(j-d+1)]^2}{ \| \phi (j-d+1) \|^2} \geq \xi \} ;
\]
clearly
$\{ j \in \Z : \;\; j \geq t_0 \} = S_{good} \cup S_{bad}  $.
Observe that this partition clearly depends on $\theta_0$, $\theta_{ab}^*$, 
etc.
We will apply Proposition 4 to analyse the closed-loop system
behaviour on $S_{good}$;
on the other hand, we will easily obtain bounds on the system behaviour on
$S_{bad}$.
Before doing so, following \cite{ccta} and \cite{rw}, we partition the time index
$\{ j \in \Z: j \geq t_0 \}$ into
intervals which oscillate between $S_{good}$ and $S_{bad}$.
To this end, it is easy to see that we can define
a (possibly infinite) sequence of intervals of the form
$[ k_i , k_{i+1} )$ satisfying:
{(i)}
$k_1 = t_0$;
{(ii)}
$[ k_i , k_{i+1} )$ either belongs to $S_{good}$ or $S_{bad}$; and
{(iii)}
if $k_{i+1} \neq \infty$ and $[ k_i , k_{i+1} )$ belongs to $S_{good}$
(respectively, $S_{bad}$), then the interval
$[ k_{i+1} , k_{i+2} )$ must belong to $S_{bad}$ (respectively,
$S_{good}$).

Now we analyse the behaviour during each interval.

\noindent
{\bf Sub-Case 2.1:}
$[ k_i , k_{i+1} )$ lies in $S_{bad}$.

Let $j \in [ k_i , k_{i+1} )$ be arbitrary.
In this case
either
$\phi (j-d+1) = 0$
or
$\frac{[ \bar{w} (j-d+1)]^2}{ \| \phi (j-d+1) \|^2} \geq \xi$
holds.
In either case we have
\be
\| \phi (j-d+1) \| \leq \frac{1}{\xi^{1/2}}  | \bar{w}(j-d+1) |, \; j \in [ k_i , k_{i+1} ) .
\label{phibd1}
\ee
From (\ref{crude}) and Proposition 2 we have
\[
\| \phi (j-d+2) \|  \leq  
\frac{c_1}{\xi^{1/2}} | \bar{w}(j-d+1) | + 
\]
\be
c_1 | y^* ( j+d+1 ) | +
c_1 | w ( j+1 ) | ,
 \; j \in [ k_i , k_{i+1} ) .
\label{bad_bound}
\ee
If we combine this with (\ref{phibd1})
we end up with
\be
\| \phi (j) \| \leq
\left\{ \begin{array}{ll}
 \frac{1}{\xi^{1/2}} | \tilde{w}(j) |  & \;\mbox{if } j=k_i \\
c_1 ( 1 + \frac{1}{ \xi^{1/2} } )  | \tilde{w}(j-1) |  & \;
\mbox{if } j = k_i+1,..., k_{i+1} .
\end{array}
\right.
\label{phibd}
\ee

\vspace{0.2cm}
\noindent
{\bf Sub-Case 2.2:}
$[ k_i , k_{i+1} )$ lies in $S_{good}$.

This case is much more involved than in the proof of \cite{ccta} and \cite{rw} since the
bound on $\| \Delta (t) \|$ provided by Claim 2
extends both forward and backward in time, occasionally outside $S_{good}$.
Hence, we need to handle the first $d$ and last $d$ time units separately.

To this end, first suppose that
$k_{i+1} \leq k_i + 2d$. Then by Claim 1 we see that there
exists a constant $\gamma_5$ so that
\[
\| \phi (k) \| \leq \gamma_5 \lambda^{k-k_i} \| \phi (k_i) \| +
\]
\be
\;\; 
\sum_{j=k_i}^{k-1} \gamma_5 \lambda^{k-1-j} | \tilde{w} (j) | , \;
k_i \leq k \leq k_{i+1} .
\label{easy}
\ee

Now suppose that $k_{i+1} > k_i+2d$. Define
$\kib = k_i +d $ and $ \kioneb = k_{i+1} -d $.
Let $j \in [ k_i , k_{i+1} )$ be arbitrary; then
\be
\rhoj \frac{ w(j-d+1)^2}{\| \phi (j-d+1) \|^2} < \xi .
\label{dbd}
\ee
Combining this with
Proposition 1 we have that
\be
\sum_{j=\uk}^{\bk-1} \nu (j)^2 
 \leq   8 \|{\cal S } \|^2  + 4 \xi ( \bk - \uk ) , \;
\kib \leq \uk < \bk \leq \kioneb .
\label{starstar}
\ee
From Claim 2 
there exists a constant $\gamma_6$ so that
\begin{eqnarray*}
\sum_{j=\uk}^{\bk-1} \| \Delta (j) \| & \leq &
  \gamma_6 ( \bk - \uk) ^{1/2} + \gamma_6 \xi^{1/2}  ( \bk - \uk ) , \\
 && \;\;
\kib \leq \uk <  \bk \leq \kioneb .
\end{eqnarray*}
Now we will apply Proposition 4: we set
\[
\beta_0 =  0 ,
\; \beta_1 =   \gamma_6 ,
\;
\beta_2 = \gamma_6 \xi^{1/2} , \;
\mu = \lambda.
\]
With 
$N$ chosen as in Case 1, we have that
$\underline{\delta} := \frac{\lambda}{ \gamma_1^{1/N}} - \lambda_1 > 0 $;
we need
${\beta_2} < \frac{1}{ \gamma_1} \underline{\delta} $,
%or equivalently
%\begin{eqnarray*}
%2 \gamma_1 \eps^{1/2} + \frac{2 \eps^{1/2}}{N} & < & \frac{1}{N \gamma_1} \underline{\delta} \\
%\Leftrightarrow \;  \eps < \frac{ \underline{\delta}^2}{4 \gamma_1^2 ( \gamma_1 N + 1 )^2} ,
%\end{eqnarray*}
which will certainly be the case if we set
$\xi :=   \frac{ \underline{\delta}^2}{2 \gamma_1^2  \gamma_6^2 }  $.
From Proposition 4 we see that there exists a constant $\gamma_7$ so that the state
transition matrix $\Phi (t , \tau )$ corresponding to ${A}_{nom}  +
\Delta (t)$ satisfies
\[
\| \Phi (t, \tau ) \| \leq \gamma_7 \lambda^{t - \tau } , \;
\kib \leq \tau \leq t \leq \kioneb .
\]
Hence, we see from (\ref{etabd}) and (\ref{starstar}) that
\[
| \eta (t) | \leq c_2 ( 1 + 8^{1/2} \| S \| + 4 \xi d ) \tilde{w} (t) , \;
\bar{k}_i \leq t \leq \uk_{i+1};
\]
if we
now apply this to
(\ref{starplant})
then we see that there exists a constant $\gamma_8$ so that
\[
\| \bphi (k) \| \leq \gamma_8 \lambda^{k- \bar{k}_i} \| \bar{\phi} ( \bar{k}_i) \| +
\sum_{j=\bar{k}_i}^{k-1}  \gamma_8 \lambda^{k-1-j} | \tilde{w} (j)|  
\]
for
$\kib \leq k \leq \kioneb  $.
We can use Claim 1 
to extend the bound to the rest of
$[ k_i , k_{i+1 } )$: 
there exists a constant $\gamma_9$ so that
\[
\| \bphi (k) \| \leq \gamma_9 \lambda^{k- k_i} \| \bphi (k_i) \| +
\sum_{j=k_i}^{k-1}  \gamma_9 \lambda^{k-1-j} | \tilde{w} (j)|  
\]
for
$k_i \leq k \leq k_{i+1}  $.
This completes Sub-Case 2.2.

Using an argument virtually identical to that used in the last part of the
proof of Theorem 1 of \cite{ccta} and \cite{rw}, we can glue the bounds from Sub-Case 1 and
Sub-Case 2 together; using Claim 1, we conclude that there exists a constant $\gamma_{10}$ so that
\[
\| \phi (k) \| \leq \gamma_{10} \lambda^{k- t_0} \| x_0 \| +
\sum_{j=t_0}^{k}  \gamma_{10} \lambda^{k-j}( | y^* (j) | + | w (j)|)
\]
for $k \geq t_0$.
This completes Case 2.

Now suppose that $w=0$.
From 
(\ref{epsdude}) and (\ref{edude}) we have
\[
	\varepsilon(t) = e(t) + \phi(t-d)^\top
	 \left[ \hat\theta(t-1) - \hat\theta(t-d) \right] , \; t \geq t_0 + d,
\]
so if
$\|\phi(t-d)\|\neq0$, we have
\[
	\frac{|\varepsilon(t)|}{\|\phi(t-d)\|} \leq 
	\frac{|e(t)|}{\|\phi(t-d)\|} +
	\Vert \hat\theta(t-1) - \hat\theta(t-d) \Vert .
\]
From the first estimator property of Proposition 1 we obtain
\[
 \frac{|\varepsilon(t)|}{\|\phi(t-d)\|} 
\leq
\sum_{j=0, \phi(t-d-j)\neq0 }^{d-1}
\frac{|e(t-j)|}{\|\phi(t-d-j)\|}  .
\]
By
Cauchy-Schwartz
we obtain
\begin{eqnarray*}
\frac{|\varepsilon(t)|^2}{\|\phi(t-d)\|^2} 
&\leq& 
\left(\sum_{j=0, \phi(t-d-j)\neq0}^{d-1}
\frac{|e(t-j)|}{\|\phi(t-d-j)\|}\right)^2 \\
%&\leq& 
%\left(\sum_{j=0, \phi(t-d-j)\neq0}^{d-1}
%\frac{|e(t-j)|^2}{\|\phi(t-d-j)\|^2}\right)\left(\sum_{j=0}^{d-1}1\right)\\
&\leq& 
d\sum_{j=0, \phi(t-d-j)\neq0}^{d-1}
\frac{|e(t-j)|^2}{\|\phi(t-d-j)\|^2}.
\end{eqnarray*}
Hence,
for $T > t_0 + 2d-1$:
\begin{eqnarray*}
&&\sum_{t=t_0+2d-1, \phi(t-d)\neq0}^{T} 
\frac{|\varepsilon(t)|^2}{\|\phi(t-d)\|^2} \\
&\leq& 
\sum_{t=t_0+2d-1, \phi(t-d)\neq0}^{T} \left(d\sum_{j=0, \phi(t-d-j)\neq0}^{d-1}
\frac{|e(t-j)|^2}{\|\phi(t-d-j)\|^2}\right) \\
%&=& 
%d \sum_{j=0, \phi(t-d-j)\neq0}^{d-1}\sum_{t=t_0+d, \phi(t-d)\neq0}^{T} 
%\frac{|e(t-j)|^2}{\|\phi(t-d-j)\|^2}  \\
%& \leq &
%d \sum_{j=0, \phi(t-d-j)\neq0}^{d-1}\sum_{t=t_0+d, \phi(t-d)\neq0}^{\infty} 
%\frac{|e(t-j)|^2}{\|\phi(t-d-j)\|^2} \\
&\leq& 
d^2 \sum_{t=t_0+d, \phi(t-d)\neq0}^{\infty} 
\frac{|e(t)|^2}{\|\phi(t-d)\|^2} \\
% &\leq& d^2 \sum_{t=0}^{\infty} \frac{|e(t-j)|^2}{\|\phi(t-d-j)\|^2} \\
&\leq&
 8d^2\|{\cal S}\|^2 \mbox{  (by Proposition 1).} 
\end{eqnarray*}
Since $\eps (t) = 0$ if $\phi (t-d) = 0$,
if we now apply the bound on $\phi (t)$ proven above, we conclude that
\begin{eqnarray*}
\sum_{t=t_0+2d-1}^{\infty}  \eps (t)^2 & \leq & 8 d^2 \| {\cal S} \|^2 \times
\sup_{j \geq t_0 } \| \phi (j) \|^2 
\end{eqnarray*}
\[
\leq  
8 d^2 \| {\cal S} \|^2 c^2 [  \| x_0 \|^2 + ( \frac{1}{1- \lambda})^2
\sup_{j \geq t_0 } | y^* (j) |^2 ] ,
\]
which yields the desired result.
$\Box$

\begin{remark}
The linear-like bound proven in Theorem 1 can be leveraged to prove that
parametric time-variations can be tolerated.
So suppose that the actual plant model is
\be
y(t+d) = \phi(t)^T \theta^* (t) + \bar{w} (t) , \; \phi (t_0) = \phi_0 ,
\label{newplantnn}
\ee
with $\theta^* (t) \in {\cal S}$ for all $t \in \Z$. We adopt a common
model of acceptable time-variations used in adaptive control:
with $c_0 \geq 0$ and $\epsilon >0$, we let
$s( {\cal S} , c_0,  \epsilon )$ denote the subset of
$\linf( \R^{n+m+d})$ whose elements $\theta^*$ satisfy
$\theta^* (t) \in {\cal S}$ for every $t \in \Z$ as well as
\[
\sum_{t=t_1}^{t_2-1} \| \theta^* (t+1) - \theta^* (t) \| \leq c_0 + \epsilon ( t_2 - t_1 ) , \;
t_2 > t_1
\]
for every $t_1 \in \Z$.
If we argue as in \cite{ccta} and \cite{rw}, we can show that
for every $c_0 \geq 0$, if $\epsilon$ is small enough
then the proposed controller will still provide
linear-like bounds on $\phi (t)$ for all $\theta^* \in 
s( {\cal S} , c_0,  \epsilon )$.
\end{remark}

\begin{remark}
The linear-like bounds proven in Theorem 1 can be used in conjunction
with the Small Gain Theorem to
prove that the
closed-loop system tolerates a degree of unmodelled dynamics.
\end{remark}

\section{A Simulation Example}

\begin{figure}
\center\includegraphics[width=0.9\columnwidth]{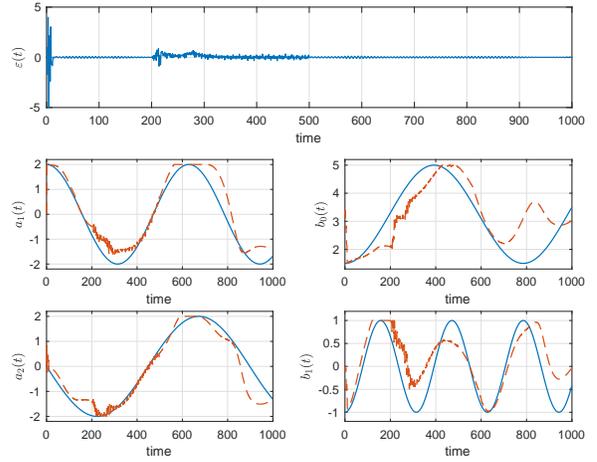}
\caption{The upper plot shows the tracking error; the lower four plots 
show the parameter estimates (dashed) as well as the actual 
parameters (solid).}
\label{fig1}
\end{figure}

Here we provide a simulation example to illustrate the results of this paper. Consider the time-varying plant
\begin{multline*} 
y(t+1)=-a_1(t)y(t)-a_2(t)y(t-1)\\+b_0(t)u(t)+b_1(t)u(t-1)+w(t)
\end{multline*}
with $a_1(t)\in[-2,2], a_2(t)\in[-2,2], b_0(t)\in[1.5,5]$ and $b_1(t)\in[-1,1]$. 
We apply the proposed adaptive controller (with $\delta = \infty$) 
to this plant when
\begin{eqnarray*}
a_1(t) =2\cos(.01t), &  a_2(t) =-2\sin(.007t), \\ 
b_0(t) =3.25-1.75\cos(.008t), & b_1(t) = - \cos(.02t),
\end{eqnarray*}
$$y^*(t)=\cos(t), \; w(t)=\left\{\begin{matrix} 
0.1 \cos(10t) && 200< t\leq 500 \\
0 && \mbox{otherwise;}
\end{matrix}\right.   $$
we set $y(-1)=y(0)=-1$, $u(-1)=0$, and the 
initial parameter estimates to the midpoint of the respective intervals. 
Figure \ref{fig1} shows the results. 
As expected, the controller does a good job of tracking when there is
no disturbance;
while the tracking degrades when the disturbance is turned on at $t=200$,
it quickly
improves when the disturbance returns to zero at $t=500$.
Furthermore, the estimator tracks the time-varying parameters fairly well.
%\begin{figure}
%\center\includegraphics[width=0.9\columnwidth]{exA_cdc2_v2}
%\caption{The upper plot shows the tracking error; the lower plot shows control input.}
%\label{fig1}
%\end{figure}
%\begin{figure}
%\center\includegraphics[width=0.9\columnwidth]{exB_cdc2_v2}
%\caption{The plots show the parameter estimates $\hat\theta(k)$ (dashed) and actual parameters $\theta^*(k)$ (solid).}
%\label{fig2}
%\end{figure}

\section{Summary and Conclusions}

Under suitable assumptions, here
we show that
if the original, ideal, projection algorithm is used in the
estimation process, then
% (subject to the common assumption that the plant 
%parameters lie in a convex,
%compact set and that the parameter estimates are restricted to that set), then  
the corresponding
$d$-step-ahead adaptive controller
{\bf guarantees linear-like convolution bounds on the closed
loop behaviour}; this confers exponential stability and a bounded
noise gain, unlike almost all other parameter adaptive controllers.
This can be leveraged in a modular way to prove tolerance to unmodelled dynamics and
plant parameter variation.

In the case of a zero disturbance, it is proven that asymptotic
tracking is achieved; we are presently working
on obtaining a
bound on the tracking quality 
in terms of the size of the disturbance.
In this approach we assumed that the sign of the high frequency 
gain is known; we are presently trying to use a multi-estimator
approach to remove this assumption.

%\section{Acknowledgements}

%This research was supported by a Grant from the Natural Sciences and
%Engineering Research Council of Canada.

\section{Appendix}

\noindent
{\bf Proof of Proposition 3:}

To proceed, we analyse the good closed-loop model of Section III.A.
From (\ref{good_model}), it is clear that we need to obtain a bound on
the terms $B_1 \eps (t+1)$ and 
$B_2 \eps ( t+j)$ for $j=1,...,d+1$.
It will be convenient to define an intermediate quantity
\footnote{It is similar to $\nu (t-1)$ except for the $\eps(t)$ rather
than $e(t)$ at the end.}:
\[
\bar{\nu} (t-1) :=
\rho_{\delta} ( \phi (t-d) , e(t)) \times
\frac{ \phi (t-d)}{\| \phi (t-d) \|^2} \eps(t) 
\]

\noindent
{\bf Step 1: Obtain a desireable bound on 
$B_i \eps (t)$ in terms of $\bar{\nu } (t-1)$, $\phi (t-d)$ and
$\bar{w} (t-d)$.}

First of all, for $i=1,2$ define
\[
\bar{\Delta}_i (t) := 
\rho_{\delta} ( \phi (t-d) , e(t)) 
\frac{ \eps (t) }{ \| \phi (t-d)\|^2} B_i \phi (t-d)^T .
\]
It is easy to see that
\[
 \bar{\Delta}_i (t) \phi (t-d) = 
\rho_{\delta} ( \phi (t-d) , e(t)) B_i { \eps (t) } .
\]
So
\begin{eqnarray}
B_i \eps (t) &=& \rho_{\delta} ( \phi (t-d) , e(t)) B_i \eps (t) + \\
&& \underbrace{[ 1 - \rho_{\delta} ( \phi (t-d) , e(t)) ] \eps (t)}_{=: \eta_0 (t)} B_i \nonumber \\
&=& \bar{\Delta}_i (t) \phi (t-d) + B_i \eta_0 (t) .
\label{eps0bd}
\end{eqnarray}
Using 
(\ref{epsdude}) and (\ref{edude})
and the definition of $\rho_{\delta} ( \phi (t-d) , e(t))$
it is easy to show that
\be
| \eta_0 (t) | \leq
\frac{ 4 \| S \| + \delta}{ \delta} | 
\bar{w} (t-d) | ;
\label{etabdd}
\ee
furthermore, it is clear that
\be
\| \bar{\Delta}_i (t) \| =
\rho_{\delta} ( \phi (t-d) , e(t)) \frac{ | \eps (t) | }{ \| \phi (t-d)\|} 
= | \bar{\nu} (t-1) | .
\label{d0bd}
\ee

\vspace{0.3cm}
\noindent
{\bf Step 2:
Bound $\bar{\nu} (t-1)$ in terms of
$ \nu (t) ,..., \nu (t-d)$.}

It follows from the formulas for $\eps (t)$ and $e(t)$ given in 
(\ref{epsdude}) and (\ref{edude}) that
\begin{eqnarray*}
\eps (t) &=& e(t) + \phi (t-d)^T [ \hat{\theta} (t-1) - \hat{\theta} (t-d) ] .
\end{eqnarray*}
Using the definition of $\bar{\nu} (t-1)$ we have
\begin{eqnarray*}
| \bar{\nu} (t-1) | & \leq & \rho_{\delta} ( \phi (t-d) , e(t)) \frac{ | e (t) | }{ \| \phi (t-d)\|} + \\
&&  \| \hat{\theta} (t-1) - \hat{\theta} (t-d) \| \\
&=& | \nu (t-1) | + \| \hat{\theta} (t-1) - \hat{\theta} (t-d) \| .
\end{eqnarray*}
Now it follows from the estimator update law that
\begin{eqnarray*}
\| \hat{\theta} (t-1) - \hat{\theta} (t-d)  \| & \leq &
\| \hat{\theta} (t-1) - \hat{\theta} (t-2 ) \| + \cdots   \\
&&
+ \| \hat{\theta} (t-d+1) - \hat{\theta} (t-d ) \|
\end{eqnarray*}
\[
    \leq 
\rho_{\delta} (  \phi (t-d-1), e (t-1) ) \frac{| e( t-1)|}
{\|  \phi (t-d-1) \| } + \cdots + 
\]
\[
\rho_{\delta} ( \phi ( t-2d+1) , e (t-d+1)) \frac{| e( t-d+1)|}
{\|  \phi (t-2d+1) \| }  
\]
\[
 \leq  \sum_{j=2}^{d} | \nu ( t-j ) |.
\]
We conclude that
\begin{eqnarray}
| \bar{\nu} (t-1) | & \leq &
\sum_{j=1}^{d} | \nu ( t-j ) |.
\label{nubbd}
\end{eqnarray}

\noindent
{\bf Step 3: Obtain a bound on $B_i \eps (t)$ in terms of
$\{ \nu (t) ,..., \nu (t-d) \}$, $\phi (t-d)$ and $\bar{w} (t-d)$.}

If we combine 
(\ref{eps0bd}), (\ref{etabdd}), (\ref{d0bd}) and (\ref{nubbd}), we see that
\[
B_i \eps (t) = \bar{\Delta}_i (t) \phi (t-d) + B_i \eta_0 (t)
\]
with
\[
\| \bar{\Delta} _i (t) \| \leq 
\sum_{j=0}^{d-1} | \nu ( t-j ) |
\]
and
\[
| \eta_0 (t) | \leq  \frac{4 \| S \| + \delta}{ \delta} | \bar{w} (t-d) | .
\]

\noindent
{\bf Step 4: Apply the result of Step 3 to (\ref{good_model}).}

We can apply the above result to each of the terms on the RHS of
(\ref{good_model}) containing $\eps ( \cdot )$.
So from Step 3 we have
\be
B_1 \eps (t+1) = \bar{ \Delta}_1 ( t+1 ) \phi (t+1-d) + B_1 \eta_0 (t+1)
\label{term1}
\ee
and
\[
\frac{a_{d-j}}{b_0} B_2 \eps ( t+1+j) =
\frac{a_{d-j}}{b_0} \bar{ \Delta}_2 ( t+1+j ) \phi (t+1-d+j) +
\]
\be
\frac{a_{d-j}}{b_0}  B_2 \eta_0 (t+1-j) , \; j =0,1,..., d.
\label{term2}
\ee
Each term except one is of the desired form: the case of $j=d$ is
problemmatic since it contains a $\phi (t+1)$ term. However, we can use
the
crude model given in (\ref{crude}) to see that
\[
\frac{a_0}{b_0} \bar{ \Delta}_2 ( t+d+1 ) \phi (t+1) =
\frac{a_0}{b_0} \bar{ \Delta}_2 ( t+d+1 ) \times
\]
\be
[ A_b (t) \phi (t) + B_3 (t) y^* (t+d+1) + B_4 w (t+1) ] .
\label{term3}
\ee
If we now combine (\ref{term1}), (\ref{term2}) and (\ref{term3}), we
see that we should define
\[
\Delta_0 (t) = \frac{a_0}{b_0} \bar{ \Delta}_2 ( t+d+1 ) A_b (t)  +
\frac{a_1}{b_0 } \bar{ \Delta}_2 ( t+d ),
\]
\[
\Delta_j (t) = \frac{a_{j+1}}{b_0} \bar{\Delta}_2 ( t+d -j ) , \;
j = 1,..., d-2 ,
\]
and
\[
\Delta_{d-1} (t) = \bar{\Delta}_1 (t+1) + \frac{a_d}{b_0} \bar{\Delta}_2 ( t +1 ).
\]
It is clear from Step 3 that this choice of $\Delta_i$ has the
desired property. 
Last of all, we group all of the remaining terms into $\eta (t)$:
\begin{eqnarray*}
\eta (t) &:=&
B_1 y^*(t+1) + B_2 \frac{1}{b_0}
\sum_{j=0}^d a_{d-j} y^* ( t+1 +j ) - \\
&& \frac{1}{b_0} B_2 w(t+d+1) + B_1 \eta_0 (t+1) + \\
&& \frac{a_0}{b_0} \bar{\Delta}_2 (t+d+1) [ B_3 (t) y^* (t+d+1) + \\
&&
B_4 w (t+1) ] + 
 B_2 \sum_{j=0}^{d-1} \frac{a_{d-j}}{b_0} \eta_0 (t+1-j).
\end{eqnarray*}
If we apply Proposition 2 and use the bound on $\eta_0 (t)$ given in Step 3,
we see that $\eta (t)$ has the desired property.

\noindent
$\Box$

\end{document}